\documentclass[11pt, A4 paper]{article}

\usepackage{amsmath}
\usepackage{times,epsfig,makeidx,amsfonts,amsmath,dcolumn,enumerate,multirow,graphicx,amssymb,colortbl,bm}
\usepackage{times}
\usepackage{natbib}

\usepackage[margin=1in]{geometry}

\newtheorem{theorem}{Theorem}

\newtheorem{corollary}{Corollary}

\newtheorem{remark}{Remark}

\newenvironment{proof}[1][Proof]{\noindent\textbf{#1.} }{\ \rule{0.5em}{0.5em}}

\begin{document}

\title{On the propriety of the posterior of hierarchical linear mixed models with flexible random effects distributions}
\author{  {\large F. J. Rubio}\footnote{{\sc University of Warwick, Department of Statistics, Coventry, CV4 7AL, UK.} E-mail:  Francisco.Rubio@warwick.ac.uk} \\
}
\date{}

\maketitle

\begin{abstract}
The use of improper priors in the context of Bayesian hierarchical linear mixed models has been studied under the assumption of normality of the random effects. We study the propriety of the posterior under more flexible distributional assumptions and general improper prior structures.
\\
\\
\noindent {\it keywords:} {Half-Cauchy prior; Improper priors; Skew--normal; Two--piece normal.}

\end{abstract}

\section{Introduction}

Hierarchical linear mixed models (LMM) are often used to account for parameter variation across groups of observations. The general formulation of this type of models is given by the model equation:

\begin{eqnarray}\label{LMM}
{\bf y} = {\bf X}\bm{\beta} + {\bf Zu} + \bm{\varepsilon},
\end{eqnarray}

\noindent where ${\bf y}$ is an $n\times 1$ vector of data, $\bm{\beta}$ is a $p\times 1$ vector of \emph{fixed effects}, ${\bf u}$ is a $q\times 1$ vector of \emph{random effects}, ${\bf X}$ and ${\bf Z}$ are known design matrices of dimension $n\times p$ and $n\times q$, respectively, and $\bm{\varepsilon}$ is an $n\times 1$ vector of residual errors. We assume that $\operatorname{rank}({\bf X})=p$ and $n>p$ throughout. The random vectors $\bm{\varepsilon}$ and ${\bf u}$ are typically assumed to be normally distributed. Given that the normality assumption can be restrictive in practice, alternative distributional assumptions have been explored, such as finite mixtures of normal distributions \citep{ZD01}, scale mixtures of skew-normal distributions \citep{L10}, and Bayesian nonparametric approaches \citep{D10}. However, in a Bayesian context with improper priors, only the model with normal assumptions has been studied. In this line, \cite{HC96} analysed the following hierarchical structure:

\begin{eqnarray}\label{Hierarchy1}
\bm{\varepsilon}\vert \sigma_0 &\sim& N_n({\bf 0},\sigma_0^2{\bf I}_n)\notag\\
{\bf y} \vert {\bf u},\bm{\beta}, \sigma_0 &\sim& N_n({\bf X}\bm{\beta} + {\bf Z}{\bf u},\sigma_0^2{\bf I}_n),\notag\\
{\bf u}\vert \sigma_1,\dots,\sigma_r &\sim& N_q({\bf 0},{\bf D}),
\end{eqnarray}

\noindent with an improper prior structure as follows:

\begin{eqnarray}\label{HCPrior}
\pi(\bm{\beta},\sigma_0,\sigma_1,\dots,\sigma_r) &\propto& \prod_{i=0}^r \dfrac{1}{\sigma_i^{2a_i+1}},
\end{eqnarray}

\noindent where $\bm{u}=({\bf u}_1^{\top},\dots,{\bf u}_r^{\top})^{\top}$, ${\bf u}_i$ is $q_i\times 1$, $\sum_{i=1}^r q_i =q$, ${\bf D}=\oplus_{i=1}^r \sigma_i^2 {\bf I}_{q_i}$, $i=1,\dots,r$, and $N_d({\bf m},{\bf S})$ denotes a $d-$variate normal distribution with mean ${\bm m}$ and covariance matrix ${\bf S}$. The $r$ subvectors of ${\bf u}$ correspond to the $r$ different random factors in the experiment. These assumptions imply that the random effects are assumed to be independent. This independence assumption restricts the applicability of this sort of models, however, it covers some models of interest in meta-analysis \citep{HC96,TH09}. \cite{HC96} obtained necessary and sufficient conditions for the propriety of the corresponding posterior distribution. This kind of improper priors may be of interest in the lack of strong prior information, although their use often represents a controversial topic. In particular, the choice $a_0=0$ and $a_1=\dots=a_r=-1/2$ is referred to as the \emph{standard diffuse prior} \citep{TH09}. \cite{S01} studied an extended improper prior structure which allows for assigning either proper or improper priors on the parameters of the distributions of the residual errors and the random effects as follows:

\begin{eqnarray}\label{SPrior}
\pi(\bm{\beta},\sigma_0,\sigma_1,\dots,\sigma_r) &\propto& \prod_{i=0}^r \dfrac{1}{\sigma_i^{2a_i+1}}\exp(-b_i/\sigma_i^2).
\end{eqnarray}

This prior corresponds to assigning an inverse gamma distribution to the $\sigma_i^2$'s when $a_i>0$ and $b_i>0$, and it contains the prior structure (\ref{Hierarchy1}) for $b_0=\dots=b_r=0$ (note that here we are presenting the priors on $(\sigma_0,\sigma_1,\dots,\sigma_r)$ rather than the equivalent priors on $(\sigma_0^2,\sigma_1^2,\dots,\sigma_r^2)$ presented in \citealp{HC96} and \citealp{S01}).

In this paper, we explore extensions of the hierarchical structure (\ref{Hierarchy1})--(\ref{SPrior}) by considering more flexible random effect distributions. In Section \ref{TPNRE}, we present an extension of this model for the case when the random effects are distributed according to a two--piece normal distribution. We show that the corresponding posterior is proper essentially under the conditions presented in \cite{S01} for the model with normal assumptions (\ref{Hierarchy1})--(\ref{SPrior}). In Section \ref{GenExtensions}, we characterise a rich class of parametric distributions, which are obtained by adding a shape parameter to the normal distribution, that also preserve the existence of the posterior distribution when they are used for modelling the random effects in (\ref{LMM}). In Section \ref{DiffusePrior}, we introduce an alternative improper prior structure with proper heavy-tailed priors on the scale parameters $(\sigma_1,\dots,\sigma_r)$. We conclude with a discussion on the implementation of this sort of models as well as possible extensions of this work.


\section{Two--piece normal random effects}\label{TPNRE}

A random variable $U\in{\mathbb R}$ is said to be distributed according to a two-piece normal distribution, denoted $U\sim \text{TPN}(\mu,\sigma,\gamma)$, if its density function can be written as \citep{AGQ05}:

\begin{eqnarray}\label{TPN}
s(u\vert\mu,\sigma,\gamma) = \dfrac{2}{\sigma[a(\gamma)+b(\gamma)]}\left[\phi\left(\dfrac{u-\mu}{\sigma b(\gamma)}\right)I(u<\mu) + \phi\left(\dfrac{u-\mu}{\sigma a(\gamma)}\right)I(u\geq\mu) \right],
\end{eqnarray}

\noindent where $\phi$ denotes the standard normal density function, $I$ denotes the indicator function, and $\{a(\cdot),b(\cdot)\}$ are positive differentiable functions. This density is continuous, unimodal, with mode at $\mu\in{\mathbb R}$, scale parameter $\sigma\in{\mathbb R}_+$, and skewness parameter $\gamma\in\Gamma\subset\mathbb R$. It coincides with the normal density when $a(\gamma)=b(\gamma)$, and it is asymmetric for $a(\gamma)\neq b(\gamma)$ while keeping the normal tails in each direction. The most common choices for $a(\cdot)$ and $b(\cdot)$ correspond to the \emph{inverse scale factors} parameterisation $\{a(\gamma),b(\gamma)\}=\{\gamma,1/\gamma\}$, $\gamma \in {\mathbb R}_+$ \citep{FS98}, and the $\epsilon-$\emph{skew} parameterisation $\{a(\gamma),b(\gamma)\}=\{1-\gamma,1+\gamma\}$, $\gamma \in (-1,1)$ \citep{MH00}.

Consider now the following extension of the hierarchical structure (\ref{Hierarchy1}),

\begin{eqnarray}\label{HierarchyTPN}
\bm{\varepsilon}\vert \sigma_0 &\sim& N_n({\bf 0},\sigma_0^2{\bf I}_n)\notag\\
{\bf y} \vert {\bf u},\bm{\beta}, \sigma_0 &\sim& N_n({\bf X}\bm{\beta} + {\bf Z}{\bf u},\sigma_0^2{\bf I}_n),\notag\\
u_{ik}\vert \sigma_i,\gamma_i &\stackrel{ind.}{\sim}& \text{TPN}(0,\sigma_i,\gamma_i ),\,\,\, k=1,\dots,q_i, \,\,\, i=1,\dots,r,
\end{eqnarray}

\noindent with prior structure

\begin{eqnarray}\label{PriorH1Ext}
\pi(\bm{\beta},\sigma_0,\sigma_1,\dots,\sigma_r,\gamma_1,\dots,\gamma_r) &\propto& \dfrac{1}{\sigma_0^{2a_0+1}}\exp(-b_0/\sigma_0^2) \prod_{i=1}^r\dfrac{\pi(\gamma_i)}{\sigma_i^{2a_i+1}}\exp(-b_i/\sigma_i^2),
\end{eqnarray}

\noindent where $\pi(\gamma_i)$ are proper priors, and ${\bf u}_i=(u_{i1},\dots,u_{iq_i})^{\top}$. The following result presents necessary and sufficient conditions for the propriety of the corresponding posterior.

\begin{theorem}\label{STP}
Consider the LMM (\ref{LMM}) with hierarchical structure (\ref{HierarchyTPN})-(\ref{PriorH1Ext}), and denote \\ $t = \operatorname{rank}\{(\bf{I-X(X^{\top}X)^{-1}X^{\top})Z}\}$. Suppose that there exist constants $m,M>0$ such that $h(\gamma)=\min\{a(\gamma),b(\gamma)\}<m$ and  $H(\gamma)=a(\gamma)+b(\gamma)>M$. Consider also the following conditions:

\begin{enumerate}[(a)]
\item[(a)] For $i=1,\dots,r$ either $a_i<  b_i=0$ or $b_i>0$,
\item[(b1)] $q_i+2a_i>0$,
\item[(b2)] $q_i+2a_i>q-t$ for all $i=1,\dots,r$,
\item[(c1)] $n - p + 2\sum_{i=0}^r a_i >0$,
\item[(c2)] $n- p + 2 a_0 + 2\sum_{i=1}^r \min(0,a_i)>0$,
\item[(d)] $\prod_{i=1}^{r} \int_{\Gamma}  \dfrac{h(\gamma_i)^{q_i+2a_i}}{[a(\gamma_i)+b(\gamma_i)]^{q_i}}\pi(\gamma_i) d\gamma_i < \infty$,
\item[(e)]  $\prod_{i=1}^{r}\int_{\Gamma}  H(\gamma_i)^{2a_i}\pi(\gamma_i)  d\gamma_i < \infty$.
\end{enumerate}

There are two cases:

\begin{enumerate}
\item If $t=q$ or if $r=1$, then the conditions (a), (b2), (c1), and (d) are necessary, while (a), (b2), (c2), and (e) are sufficient for the propriety of the posterior for almost any sample.

\item If $t<q$ and $r>1$, then the conditions (a), (b1), (c1), and (d) are necessary, while (a), (b2), (c2), and (e) are sufficient for the propriety of the posterior for almost any sample.

\end{enumerate}

\proof See appendix.

\end{theorem}

Conditions (a)--(c2) are exactly the same as those presented in Theorem 2 from \cite{S01} for the propriety of the posterior associated to the hierarchical structure (\ref{Hierarchy1})--(\ref{SPrior}), while the conditions (d)--(e) are specific to the skewness parameter of the two-piece normal extension proposed here. Conditions (a)--(c2) refer to the need for the presence of certain amount of  repeated measurements, while conditions (d)--(e) are trivially satisfied for some parameterisations, such as the $\epsilon-$skew parameterisation where $H(\gamma)\equiv 2$. The $\epsilon-$skew and the inverse scale factors parameterisations satisfy the boundedness conditions in Theorem \ref{STP}. Note also that the bounds $m$ and $M$ do not appear in conditions (a)--(e) but they are essential for the proof of the result. The function $H(\gamma)$ in condition (e) can be replaced by $H(\gamma)=\max\{a(\gamma),b(\gamma)\}$, which might facilitate checking this condition under some parameterisations.

A discussion on the meaning of the phrase ``for almost any sample'' in Theorem \ref{STP} seems appropriate. The propriety results presented in Theorem 2 from \cite{S01} are satisfied for the case when $2b_0 + \text{SSE}>0$, where $\text{SSE}$ represents the sum of square errors of model (\ref{LMM}). This condition, as discussed by \cite{S01}, is satisfied with probability one (this is, for almost any sample) given that the distributions involved in this model are continuous. Given that the proof of Theorem \ref{STP} is based on that of Theorem 2 from \cite{S01}, this (mild) condition is inherited by our proposed extension.

The elicitation of the proper priors $\pi(\gamma)$ can be made following the strategy proposed in \cite{RS14}. They propose assigning a uniform prior on a measure of skewness which is an injective function of the parameter $\gamma$. This induces a proper prior on $\gamma$ which reflects vague prior information on a quantity of interest. For the $\epsilon-$skew parameterisation this strategy simply leads to a uniform prior on $\gamma\in(-1,1)$.

\section{Sufficient conditions for other extensions}\label{GenExtensions}

A natural question is whether or not there are other distributions with shape parameters that can be used for modelling the random effects in (\ref{LMM}) while preserving the existence of the posterior distribution. In order to answer this question, we consider a general extension of the normal distribution obtained by using the representation proposed in \cite{FS06}. \cite{FS06} show that the transformed density function

\begin{eqnarray}\label{GenRep}
s(u\vert\mu,\sigma,\lambda) = \dfrac{1}{\sigma}p\left[\Phi\left(\dfrac{u-\mu}{\sigma}\right)\Bigg\vert\lambda\right]\phi\left(\dfrac{u- \mu}{\sigma}\right),
\end{eqnarray}

\noindent where $p:[0,1]\rightarrow{\mathbb R}_+$ is a continuous density function with shape parameter $\lambda\in\Lambda\subseteq{\mathbb R}$, and $\Phi$ denotes the standard normal distribution function, is also a density function. Moreover, for any continuous density $s$ with support on ${\mathbb R}$, there exists a density function $p$ such that the relation (\ref{GenRep}) holds. The density $s$ in (\ref{GenRep}) is symmetric if and only if $p$ is symmetric around $1/2$, and it coincides with the normal density when $p$ is the uniform density function. Throughout we denote by $U\sim\text{FSN}(\mu,\sigma,\lambda;p)$ the representation (\ref{GenRep}) of the distribution of an absolutely continuous random variable $U$ (FSN is used here as an acronym of the {\bf F}erreira-{\bf S}teel representation of transformations of the {\bf N}ormal distribution). The following result provides sufficient conditions on $p$ for the existence of the posterior of model (\ref{LMM}) with an improper prior structure and non-normal assumptions on the random effects.

\begin{theorem}\label{SE}
Consider the LMM (\ref{LMM}) with the following hierarchical structure:

\begin{eqnarray}\label{HierarchyGen}
\bm{\varepsilon}\vert \sigma_0 &\sim& N_n({\bf 0},\sigma_0^2{\bf I}_n)\notag\\
{\bf y} \vert {\bf u},\bm{\beta}, \sigma_0 &\sim& N_n({\bf X}\bm{\beta} + {\bf Z}{\bf u},\sigma_0^2{\bf I}_n),\notag\\
{u}_{ik}\vert \sigma_i,\lambda_i &\stackrel{ind.}{\sim}& \text{FSN}(0,\sigma_i,\lambda_i;p),\,\,\, k=1,\dots,q_i, \,\,\, i=1,\dots,r,
\end{eqnarray}

\noindent with prior structure:

\begin{eqnarray}\label{PriorGen}
\pi(\bm{\beta},\sigma_0,\sigma_1,\dots,\sigma_r,\lambda_1,\dots,\lambda_r) &\propto& \dfrac{1}{\sigma_0^{2a_0+1}}\exp(-b_0/\sigma_0^2) \prod_{i=1}^r \dfrac{\pi(\lambda_i)}{\sigma_i^{2a_i+1}} \exp(-b_i/\sigma_i^2).
\end{eqnarray}

Suppose that the priors $\pi(\lambda_i)$, $i=1,\dots,r$, are proper and that the density $p$ in (\ref{GenRep}) is upper bounded. Then, the conditions (a), (b2), and (c2) in Theorem \ref{STP} are sufficient for the propriety of the corresponding posterior distribution for almost any sample.

\proof See appendix.

\end{theorem}

This result implies that the use of skew-normal random effects (\citealp{A85}, whose representation is obtained for the choice $p[\Phi(u)\vert\lambda]=2\Phi(\lambda u)$, $\lambda\in{\mathbb R}$) in the hierarchical structure (\ref{HierarchyGen}) leads to a proper posterior. This result also includes other symmetric (non-normal) distributions for a symmetric (with respect to 1/2) and bounded choice of $p$. Other distributions for which the corresponding $p$ is bounded are discussed in \cite{MR13}. The density $p$ associated to the representation (\ref{GenRep}) of the two-piece normal density (\ref{TPN}) is not upper bounded \citep{MR13}. This indicates that Theorem \ref{STP} is not a consequence of Theorem \ref{SE}, which justifies their separate study. It also implies that the boundedness of the density $p$ is not a necessary condition for the existence of the posterior.

\section{An alternative weakly informative prior structure}\label{DiffusePrior}

The use of improper priors for the scale parameters in hierarchical models is certainly a debatable topic. However, the use of certain proper priors in this context has also foster some discussions. In particular, the inverse gamma distribution, which appears in the structure (\ref{SPrior}), has been shown to be influential on the shape of the posterior distribution \citep{G06}. This is an undesirable property if one is interested on using weakly informative priors (which is usually the aim when using jointly improper priors). An alternative proper prior for scale parameters that has gained a lot of attention is the half-Cauchy prior. \cite{PS12} show that this prior induces a posterior distribution with good frequentist properties. The following result presents conditions for the propriety of the posterior distribution induced by an improper prior structure with half-Cauchy priors on the scale parameters.

\begin{corollary}\label{VagueTP}
Consider the LMM (\ref{LMM}) with hierarchical structure (\ref{HierarchyTPN}) and prior structure

\begin{eqnarray}\label{PriorVagueTPN}
\pi(\bm{\beta},\sigma_0,\sigma_1,\dots,\sigma_r,\gamma_1,\dots,\gamma_r) &\propto& \dfrac{1}{\sigma_0^{2a_0+1}} \prod_{i=1}^r\dfrac{\pi(\gamma_i)}{1+\frac{\sigma_i^{2}}{s_i^2}},
\end{eqnarray}

\noindent where $\pi(\gamma_i)$ are proper priors, and $s_i>0$, $i=1,\dots,r$. Let $(X:Z)$ denote the entire design matrix. Then, the following conditions are sufficient, for almost any sample, for the propriety of the posterior distribution under any parameterisation $\{a(\gamma),b(\gamma)\}$:
\begin{enumerate}[(a)]
\item $a_0\geq 0$,
\item $\operatorname{rank}(X:Z)<n$.
\end{enumerate}

\proof See appendix.
\end{corollary}

It is worth emphasising that if we use a Dirac delta prior for $\gamma_i$ concentrated at a value $\gamma_o$ such that $a(\gamma_o)=b(\gamma_o)$, we obtain a hierarchical model with normal residual errors and normal random effects with an alternative prior structure to that in (\ref{SPrior}). This new prior structure assigns half-Cauchy priors on the scale parameters $(\sigma_1,\dots,\sigma_r)$ which, as mentioned before, has been argued to be less influential on the posterior distribution than the inverse gamma prior in (\ref{SPrior}). This feature may be appealing to practitioners in contexts with little prior information.

The prior structure in (\ref{PriorVagueTPN}) involves half-Cauchy priors with arbitrary positive scale parameters $s_i$. This can be useful to conduct sensitivity analyses on the choice of these hyperparameters.

Theorem 1 from \cite{FOS97} can be used to construct other flexible hierarchical models as follows. By using this result, we can assign an \emph{arbitrary} random effects distribution ${\bf u}\vert \bm{\theta} \sim F(\cdot\vert \bm{\theta})$, where $\bm{\theta}$ represents the parameter vector of the distribution $F$, for the LMM (\ref{LMM}) with normal residual errors $\bm{\varepsilon}\vert \sigma_0 \sim N_n({\bf 0},\sigma_0^2{\bf I}_n)$ and prior structure:

\begin{eqnarray*}\label{PriorVagueGen}
\pi(\bm{\beta},\sigma_0,\bm{\theta}) &\propto& \dfrac{\pi(\bm{\theta})}{\sigma_0^{2a_0+1}},
\end{eqnarray*}

\noindent where $\pi(\bm{\theta})$ is any proper prior. Conditions (a) and (b) in Corollary \ref{VagueTP} are sufficient for the propriety of the corresponding posterior distribution. This covers the class of random effects distributions studied in Section \ref{GenExtensions}, finite mixtures of normals, scale mixtures of normals, and skewed scale mixtures of normals with proper priors on the corresponding parameters.


\section{Discussion}

The main message of the paper is that non-normal distributions can also be used to model the random effects in hierarchical linear mixed models as well as other improper prior structures without the need for additional onerous conditions for the propriety of the posterior distribution. This is useful to construct models that are robust to departures of normality of the random effects.

Given that the proposed hierarchical structures (\ref{HierarchyTPN}) and (\ref{HierarchyGen}) contain the model with normal random effects as a particular case (for certain choices of $p$), the appropriateness of the normality assumption can be evaluated via Bayes factors since only the common parameters have improper priors. In particular, the use of the Savage--Dickey density ratio can be easily implemented to do so, as long as we are able to obtain a posterior sample from the corresponding parameters.

The use of non-normal distributional assumptions on the random effects in the context of LMM usually complicates their implementation. A general approach for sampling from the posterior distribution of the proposed hierarchical structures consists of using a Metropolis within Gibbs algorithm \citep{RR09}. However, this can be computationally demanding in moderate or high dimensions. The development of Gibbs samplers for the hierarchical structures (\ref{HierarchyTPN})--(\ref{HierarchyGen}), as well as convergence analyses of these, represent interesting research directions. In the context of LMM with normal random effects these ideas have been studied in \cite{HC96} and \cite{RH12}.

Further research includes the study of LMM with flexible distributions for the residual errors.

\subsubsection*{The case of multivariate scale mixture of normals and improper priors: a warning case?}

Another type of extension that can be of interest in practice consists of using distributions with heavier tails than those of the normal distribution for modelling the random effects. A natural candidate for this choice is the family of multivariate scale mixtures of normals, which is a rich family that includes symmetric models with heavier tails than normal. Recall first that a $p-$variate, $p\geq1$, scale mixture of normal distributions (denoted $\text{SMN}_p(\bm{\mu},\bm{\Sigma},\delta;H)$) is defined as

\begin{eqnarray*}
g({\bf x}) = \int_0^{\infty} \dfrac{\tau^{p/2}}{\det{(2\pi \bm{\Sigma})^{1/2}}} \exp\left[-\dfrac{\tau}{2}({\bf x}-\bm{\mu})^{\top}\bm{\Sigma}^{-1}({\bf x}-\bm{\mu})\right]dH(\tau\vert\delta).
\end{eqnarray*}

\noindent where $\bm{\Sigma}$ is a symmetric positive definite matrix and $H$ is a mixing distribution. For example, when $H(\tau\vert\delta)$ is a Gamma distribution with parameters $(\delta/2,\delta/2)$, we obtain the $p-$variate Student-$t$ distribution with $\delta$ degrees of freedom.

Consider now the LMM (\ref{LMM}) with the following hierarchical structure:

\begin{eqnarray}\label{HierarchySMN}
\bm{\varepsilon}\vert \sigma_0 &\sim& N_n({\bf 0},\sigma_0^2{\bf I}_n)\notag\\
{\bf y} \vert {\bf u},\bm{\beta}, \sigma_0 &\sim& N_n({\bf X}\bm{\beta} + {\bf Z}{\bf u},\sigma_0^2{\bf I}_n),\notag\\
{\bf u}_i \vert \sigma_i &\stackrel{ind.}{\sim}& \text{SMN}_{q_i}({\bf 0},\sigma_i^2 {\bf I}_{q_i},\delta;H),\,\,\, i=1,\dots,r, \notag\\
\pi(\bm{\beta},\sigma_0,\sigma_1,\dots,\sigma_r,\delta) &\propto& \pi(\delta)\prod_{i=0}^r \dfrac{1}{\sigma_i^{2a_i+1}},
\end{eqnarray}

\noindent where $\pi(\delta)$ is a proper prior on the shape parameter $\delta\in \Delta \subseteq {\mathbb R}$. It is possible to prove that for this model the marginal likelihood of the data can be written as follows:

\begin{eqnarray*}
\tilde{m}({\bf y}) = m_N({\bf y}) \int_{{\mathbb R}_+^r\times\Delta} \left[\prod_{i=1}^r \tau_i^{-a_i} dH(\tau_i\vert\delta)\right] \pi(\delta) d\delta,
\end{eqnarray*}

\noindent where $m_N({\bf y})$ is the marginal likelihood of the data with normal assumptions (\ref{Hierarchy1})-(\ref{HCPrior}). This result indicates that $\tilde{m}({\bf y})$ is proportional to the marginal likelihood of the data associated to the model with normal assumptions. Moreover, the proportionality constant depends only on the choice of the hyperparameters $a_i$, the mixing distribution $H$, and the prior $\pi(\delta)$. Therefore, conducting model selection using Bayes factors between models of the type (\ref{HierarchySMN}) with different mixing distributions (\emph{e.g.~}against the model with normal random effects) is not driven by data (!). This result provides a warning on the use of multivariate scale mixtures of normals combined with improper priors in this context. Intuitively, this result indicates that the data do not contain enough information to distinguish between models of type (\ref{HierarchySMN}).

\subsubsection*{Extensions with other generalised linear mixed models}

In the context of generalised linear mixed models (GLMM) with improper priors, only models with normal random effects appear to have been studied \citep{CH02,SL10}. Some limitations induced by the assumption of normality of the random effects have been analysed in \cite{L08}. Given that the proof technique employed in Theorems \ref{STP} and \ref{SE} is based on obtaining upper and lower bounds for the random effects distributions in terms of the normal distribution, this idea can also be applied to other GLMM with improper priors. In particular, the one-way random effect probit model studied in \cite{SL10}:

\begin{eqnarray*}
y_{ik}\vert \mu, u_i &\sim & \text{Bernoulli}[\Phi(\mu + u_i)],\,\,\, k=1,\dots,q_i, \,\,\, i=1,\dots,r,\notag\\
u_i\vert \sigma &\stackrel{ind.}{\sim}& \text{N}(0,\sigma^2),\notag \\
\pi(\mu,\sigma) &\propto& \dfrac{1}{\sigma^{2a_1+1}},
\end{eqnarray*}

\noindent can be extended to the use of two-piece normal random effects as follows.

\begin{remark}\label{TPNProbitRemark}
Consider the following one-way random effect probit model,

\begin{eqnarray}\label{TPNProbit}
y_{ik}\vert \mu, u_i &\sim & \text{Bernoulli}[\Phi(\mu + u_i)],\,\,\, k=1,\dots,q_i, \,\,\, i=1,\dots,r,\notag\\
u_i\vert \sigma,\gamma &\stackrel{ind.}{\sim}& \text{TPN}(0,\sigma,\gamma),\notag \\
\pi(\mu,\sigma,\gamma) &\propto& \dfrac{\pi(\gamma)}{\sigma^{2a_1+1}},
\end{eqnarray}

\noindent where $\pi(\gamma)$ is a proper prior. The posterior distribution of $(\mu,\sigma,\gamma)$ is proper under the following conditions:

\begin{enumerate}[(i)]
\item for each $i$, $i=1,2,\dots,r_1$, ($2\leq r_1\leq r$), there is at least one success and one failure,
\item $-(r_1-1)/2<a_1<0$,
\item  $\int_{\Gamma} H(\gamma)^{2a_1}\pi(\gamma)  d\gamma< \infty$.
\end{enumerate}

\proof The result follows by using the proof technique in Theorem \ref{STP} together with Theorem 2.1 from \cite{SL10}.
\end{remark}

In a similar way, this model can be extended to the class of random effects distributions discussed in Section \ref{GenExtensions}, which, in particular, allows for the use of skew-normal random effects.

\begin{remark}\label{GenProbitRemark}
Consider the following one-way random effect probit model,

\begin{eqnarray}\label{GenNProbit}
y_{ik}\vert \mu, u_i &\sim & \text{Bernoulli}[\Phi(\mu + u_i)],\,\,\, k=1,\dots,q_i, \,\,\, i=1,\dots,r,\notag\\
u_i\vert \sigma,\lambda &\stackrel{ind.}{\sim}& \text{FSN}(0,\sigma,\lambda;p),\notag \\
\pi(\mu,\sigma,\lambda) &\propto& \dfrac{\pi(\lambda)}{\sigma^{2a_1+1}},
\end{eqnarray}

\noindent where $\pi(\lambda)$ is a proper prior and the density $p$ is upper bounded. The posterior distribution of $(\mu,\sigma,\lambda)$ is proper under conditions (i) and (ii) from Remark \ref{TPNProbitRemark}.

\proof The result follows by using the proof technique in Theorem \ref{SE} together with Theorem 2.1 from \cite{SL10}.
\end{remark}

Other types of GLMM with improper priors where the assumption of normality of the random effects can be relaxed by using the ideas in this paper are presented in Section 4 of \cite{CH02}. The corresponding extensions and propriety results are omitted for the sake of space.


\section*{Acknowledgement}
FJR gratefully acknowledges support from EPSRC grant EP/K007521/1. I thank an Associate Editor and a referee for helpful comments.

\section*{Appendix}

\subsection*{Proof of Theorem \ref{STP}}
Let us denote $\bm{\sigma}=(\sigma_1,\dots,\sigma_r)$ and $\bm{\gamma}=(\gamma_1,\dots,\gamma_r)$. The marginal density of the data is given by

\begin{eqnarray*}
m({\bf y}) &\propto& \int f({\bf y}\vert {\bf u},\bm{\beta}, \sigma_0) f({\bf u}\vert \bm{\sigma},\bm{\gamma}) \dfrac{1}{\sigma_0^{2a_0+1}}\exp(-b_0/\sigma_0^2) \prod_{i=1}^r\dfrac{\pi(\gamma_i)}{\sigma_i^{2a_i+1}}\exp(-b_i/\sigma_i^2)  d{\bf u} d\bm{\beta} d\sigma_0 d\bm{\sigma} d\bm{\gamma}\\
&\propto& \int \dfrac{1}{\sigma_0^n } \prod_{j=1}^n \phi\left(\dfrac{y_j-{\bf x}_j^{\top}\bm{\beta} - {\bf z}_j^{\top}{\bf u}}{\sigma_0}\right) \\
&\times& \prod_{i=1}^r\prod_{k=1}^{q_i} \dfrac{1}{\sigma_i} \Biggl\{\exp\left[-\dfrac{u_{ik}^2}{2\sigma_{i}^2b(\gamma_i)^2}\right] I(u_{ik}<0)  +  \exp\left[-\dfrac{u_{ik}^2}{2\sigma_{i}^2a(\gamma_i)^2}\right] I(u_{ik}\geq0)\Biggr\}\\
&\times& \dfrac{1}{\sigma_0^{2a_0+1}}\exp(-b_0/\sigma_0^2) \prod_{i=1}^r\dfrac{\pi(\gamma_i)}{\sigma_i^{2a_i+1} [a(\gamma_i)+b(\gamma_i)]^{q_i}}\exp(-b_i/\sigma_i^2)
 d{\bf u} d\bm{\beta} d\sigma_0 d\bm{\sigma}  d\bm{\gamma} .
\end{eqnarray*}

Denote $h(\gamma)=\min\{a(\gamma),b(\gamma)\}$, then we can obtain a lower bound, up to a proportionality constant, for $m({\bf y})$ as follows

\begin{eqnarray*}
&&\int \dfrac{1}{\sigma_0^n } \prod_{j=1}^n \phi\left(\dfrac{y_j-{\bf x}_j^{\top}\bm{\beta} - {\bf z}_j^{\top}{\bf u}}{\sigma_0}\right)
 \prod_{i=1}^r\prod_{k=1}^{q_i} \dfrac{1}{\sigma_i} \exp\left[-\dfrac{u_{ik}^2}{2\sigma_{i}^2h(\gamma_i)^2}\right] \\
&\times& \dfrac{1}{\sigma_0^{2a_0+1}}\exp(-b_0/\sigma_0^2) \prod_{i=1}^r\dfrac{\pi(\gamma_i)}{\sigma_i^{2a_i+1} [a(\gamma_i)+b(\gamma_i)]^{q_i}}\exp(-b_i/\sigma_i^2)
 d{\bf u} d\bm{\beta} d\sigma_0 d\bm{\sigma} d\bm{\gamma} .
\end{eqnarray*}

By using the change of variable $\vartheta_i  = \sigma_i h(\gamma_i)$ and the lower boundedness of $h(\cdot)$ we obtain the following lower bound for $m({\bf y})$:

\begin{eqnarray}\label{STPNec}
&&\int \dfrac{1}{\sigma_0^n } \prod_{j=1}^n \phi\left(\dfrac{y_j-{\bf x}_j^{\top}\bm{\beta} - {\bf z}_j^{\top}{\bf u}}{\sigma_0}\right)
 \prod_{i=1}^r\prod_{k=1}^{q_i} \dfrac{1}{\vartheta_i} \exp\left[-\dfrac{u_{ik}^2}{2\vartheta_i^2}\right] \notag\\
&\times& \dfrac{1}{\sigma_0^{2a_0+1}}\exp(-b_0/\sigma_0^2) \prod_{i=1}^{r} \dfrac{1}{\vartheta_i^{2a_i+1}}\exp(-b_i m^2/\vartheta_i^2)
 d{\bf u} d\bm{\beta} d\sigma_0 d\bm{\vartheta} \notag\\
 &\times& \int_{\Gamma^{r}} \prod_{i=1}^{r} \dfrac{h(\gamma_i)^{q_i+2a_i}}{[a(\gamma_i)+b(\gamma_i)]^{q_i}}\pi(\gamma_i)  d\bm{\gamma}.
\end{eqnarray}

Now, denote $H(\gamma)=a(\gamma)+b(\gamma)$, then we can obtain an upper bound, up to a proportionality constant, for $m({\bf y})$ as follows

\begin{eqnarray*}
&&\int \dfrac{1}{\sigma_0^n } \prod_{j=1}^n \phi\left(\dfrac{y_j-{\bf x}_j^{\top}\bm{\beta} - {\bf z}_j^{\top}{\bf u}}{\sigma_0}\right)
 \prod_{i=1}^r\prod_{k=1}^{q_i} \dfrac{1}{\sigma_i} \exp\left[-\dfrac{u_{ik}^2}{2\sigma_{i}^2H(\gamma_i)^2}\right] \\
&\times& \dfrac{1}{\sigma_0^{2a_0+1}}\exp(-b_0/\sigma_0^2) \prod_{i=1}^{r} \dfrac{\pi(\gamma_i)}{\sigma_i^{2a_i+1}H(\gamma_i)^{q_i}} \exp(-b_i/\sigma_i^2)
 d{\bf u} d\bm{\beta} d\sigma_0 d\bm{\sigma} d\bm{\gamma} .
\end{eqnarray*}

By using the change of variable $\vartheta_i = \sigma_i H(\gamma_i)$ we obtain the following upper bound

\begin{eqnarray}\label{STPSuf}
&&\int \dfrac{1}{\sigma_0^n } \prod_{j=1}^n \phi\left(\dfrac{y_j-{\bf x}_j^{\top}\bm{\beta} - {\bf z}_j^{\top}{\bf u}}{\sigma_0}\right)
 \prod_{i=1}^r\prod_{k=1}^{q_i} \dfrac{1}{\vartheta_i} \exp\left[-\dfrac{u_{ik}^2}{2\vartheta_i^2}\right] \notag\\
&\times&  \dfrac{1}{\sigma_0^{2a_0+1}}\exp(-b_0/\sigma_0^2) \prod_{i=1}^{r} \dfrac{1}{\vartheta_i^{2a_i+1}} \exp(-b_i M^2/\vartheta_i^2)
 d{\bf u} d\bm{\beta} d\sigma_0 d\bm{\vartheta}
 \times \int_{\Gamma^{r}} \prod_{i=1}^{r} H(\gamma_i)^{2a_i}\pi(\gamma_i) d\bm{\gamma}.
\end{eqnarray}

From (\ref{STPNec}) we can identify the first factor as the marginal likelihood of the data under normal residual errors $\bm{\varepsilon}$ and normal random effects ${\bf u}$. From (\ref{STPSuf}) we can also identify the first factor as the marginal likelihood of the data under normal residual errors $\bm{\varepsilon}$ and normal random effects ${\bf u}$. Then, by combining (\ref{STPNec}), (\ref{STPSuf}) and Theorem 2 from \cite{S01} we obtain necessary and sufficient conditions for the finiteness of $m({\bf y})$.

\subsection*{Proof of Theorem \ref{SE}}
Let us denote $\bm{\sigma}=(\sigma_1,\dots,\sigma_r)$ and $\bm{\lambda}=(\lambda_1,\dots,\lambda_r)$. By using that $p$ in (\ref{GenRep}) is upper bounded, the marginal density of the data can be upper bounded, up to a proportionality constant, after integrating $\bm{\lambda}$ as follows:

\begin{eqnarray}\label{UpperBoundGR}
m({\bf y}) &\propto& \int f({\bf y}\vert {\bf u},\bm{\beta}, \sigma_0) f({\bf u}\vert \bm{\sigma},\bm{\lambda}) \dfrac{1}{\sigma_0^{b+1}} \prod_{i=1}^{r} \dfrac{\pi(\lambda_i)}{\sigma_i^{a_i+1}}  d{\bf u} d\bm{\beta} d\sigma_0 d\bm{\sigma} d\bm{\lambda}\notag\\
&\propto& \int \dfrac{1}{\sigma_0^n } \prod_{j=1}^n \phi\left(\dfrac{y_j-{\bf x}_j^{\top}\bm{\beta} - {\bf z}_j^{\top}{\bf u}}{\sigma_0}\right)
\times \prod_{i=1}^r\prod_{k=1}^{q_i} s(u_{ik}\vert 0,\sigma_i,\lambda_i) \notag\\
&\times&   \dfrac{1}{\sigma_0^{2a_0+1}}\exp(-b_0/\sigma_0^2) \prod_{i=1}^r \dfrac{\pi(\lambda_i)}{\sigma_i^{2a_i+1}} \exp(-b_i/\sigma_i^2)
 d{\bf u} d\bm{\beta} d\sigma_0 d\bm{\sigma} d\bm{\lambda} \notag\\
&\stackrel{\cdot}{\leq}& \int \dfrac{1}{\sigma_0^n } \prod_{j=1}^n \phi\left(\dfrac{y_j-{\bf x}_j^{\top}\bm{\beta} - {\bf z}_j^{\top}{\bf u}}{\sigma_0}\right)  \prod_{i=1}^r\prod_{k=1}^{q_i}\dfrac{1}{\sigma_{i}} \phi\left(\dfrac{u_{ik}}{\sigma_{i}}\right)  \notag\\
&\times&  \dfrac{1}{\sigma_0^{2a_0+1}}\exp(-b_0/\sigma_0^2) \prod_{i=1}^r \dfrac{1}{\sigma_i^{2a_i+1}} \exp(-b_i/\sigma_i^2)
 d{\bf u} d\bm{\beta} d\sigma_0 d\bm{\sigma}  .
\end{eqnarray}

Therefore, $m({\bf y})$ is upper bounded, up to a proportionality constant, by the marginal likelihood associated to a model with normal residual errors $\bm{\varepsilon}$ and normal random effects ${\bf u}$. Consequently, the sufficient conditions for the finiteness of this upper bound are also sufficient for the finiteness of $m({\bf y})$, which are obtained from Theorem 2 in \cite{S01}, and the result follows.

\subsection*{Proof of Corollary \ref{VagueTP}}

The result follows by noting that the prior (\ref{PriorVagueTPN}) assigns proper priors on the parameters of the distribution of the random effects. This propriety implies that the marginal distribution of the random effects, say $p({\bf u})$, is proper. Under these conditions, \cite{FOS97} (Theorem 1) show that conditions (a) and (b) are sufficient for the existence of the corresponding posterior distribution for almost any sample. In this case, the phrase ``for almost any sample'' refers to all the vectors ${\bf y}\in{\mathbb R}^n$ such that there is no solution for the equation ${\bf y} = {\bf X}\bm{\beta} + {\bf Zu}$.


\end{document}